\documentclass[a4paper,12pt]{amsart}
\usepackage{amssymb}
\usepackage{url}
\usepackage[T1]{fontenc}

\input xypic
\xyoption{all}
\xyoption{arc}

\newtheorem{theo}{Th\'eor\`eme}[section]

\newtheorem{prop}[theo]{Proposition}

\newtheorem{lem}[theo]{Lemme}

\newtheorem{coro}[theo]{Corollaire}

\def\remark#1{{\refstepcounter{theo}\label{#1}\noindent\sc Remarque
\arabic{section}.\arabic{theo} - }}

\def\equat{\refstepcounter{theo}$$~}
\def\endequat{\leqno{\boldsymbol{(\arabic{section}.\arabic{theo})}}~$$}

    \def\FM{{\mathbb{F}}}

    \def\NM{{\mathbb{N}}}

    \def\ZM{{\mathbb{Z}}}

\def\Ab{{\mathbf A}}    
\def\Bb{{\mathbf B}}    
    
  \def\db{{\mathbf d}}

\def\Gb{{\mathbf G}}  \def\gb{{\mathbf g}}  
\def\Hb{{\mathbf H}}

\def\Lb{{\mathbf L}}    
  \def\mb{{\mathbf m}}  
  \def\nb{{\mathbf n}}

\def\Rb{{\mathbf R}}    
\def\Sb{{\mathbf S}}  \def\sb{{\mathbf s}}  
\def\Tb{{\mathbf T}}    
\def\Ub{{\mathbf U}}    
    
  \def\wb{{\mathbf w}}  
\def\Xb{{\mathbf X}}  \def\xb{{\mathbf x}}  
\def\Yb{{\mathbf Y}}  \def\yb{{\mathbf y}}  
\def\Zb{{\mathbf Z}}  \def\zb{{\mathbf z}}

    \def\BCB{{\boldsymbol{\mathcal{B}}}}

    \def\UCB{{\boldsymbol{\mathcal{U}}}}

          \def\sdo{{\dot{s}}}

          \def\wdo{{\dot{w}}}
          \def\xdo{{\dot{x}}}

\def\Sba{{\bar{S}}}

\def\Xbt{{\tilde{\Xb}}}          
\def\Ybt{{\tilde{\Yb}}}

\def\wbdo{{\dot{\wb}}}          
\def\xbdo{{\dot{\xb}}}

          \def\Xbo{{\overline{\Xb}}}

\def\a{\alpha}
\def\b{\beta}
\def\g{\gamma}

\def\D{\Delta}

\def\ph{\varphi}
\def\l{\lambda}

\def\s{\sigma}

\def\th{\theta}

\def\t{\tau}
\def\x{\xi}
\def\z{\zeta}

\def\betb{{\boldsymbol{\beta}}}         
\def\gamb{{\boldsymbol{\gamma}}}

\def\pib{{\boldsymbol{\pi}}}            \def\pit{{\tilde{\pi}}}

\def\sigb{{\boldsymbol{\sigma}}}

\def\xib{{\boldsymbol{\xi}}}

\def\pibt{{\tilde{\boldsymbol{\pi}}}}

\DeclareMathOperator{\diag}{{\mathrm{diag}}}

\DeclareMathOperator{\Id}{{\mathrm{Id}}}

\DeclareMathOperator{\Ker}{{\mathrm{Ker}}}

\def\to{\rightarrow}
\def\longto{\longrightarrow}

\def\longtrait#1{\hspace{0.2em}\stackrel{\stackrel{\SS{#1}}{\rule{0.7cm}{0.5pt}}}{~}
\hspace{0.2em}}

\def\fonction#1#2#3#4#5{\begin{array}{rccl}
{#1} : & {#2} & \longto & {#3} \\
& {#4} & \longmapsto & {#5} 
\end{array}}

\def\fonctio#1#2#3#4{\begin{array}{ccc}
{#1} & \longto & {#2} \\
{#3} & \longmapsto & {#4} 
\end{array}}

\def\vide{\varnothing}
\def\DS{\displaystyle}
\def\SS{\scriptstyle}

\def\finl{~$\SS \square$}

\def\infspe{\hspace{0.1em}\mathop{\preccurlyeq}\nolimits\hspace{0.1em}}

\def\lexp#1#2{\kern\scriptspace\vphantom{#2}^{#1}\kern-\scriptspace#2}
\def\le{\hspace{0.1em}\mathop{\leqslant}\nolimits\hspace{0.1em}}
\def\ge{\hspace{0.1em}\mathop{\geqslant}\nolimits\hspace{0.1em}}

\mathchardef\inferieur="321E
\mathchardef\superieur="321F

\def\eqna{\begin{eqnarray*}}
\def\endeqna{\end{eqnarray*}}

\def\itemth#1{\item[${\mathrm{(#1)}}$]}

\def\gfp{{\FM_{\! p}}}

\catcode`\@=11
\long\def\@car#1#2\@nil{#1}
\long\def\@first#1#2{#1}
\long\def\@second#1#2{#2}
\long\def\ifempty#1{\expandafter\ifx\@car#1@\@nil @\@empty
  \expandafter\@first\else\expandafter\@second\fi}
\catcode`\@=12

\def\matrice#1{\begin{pmatrix} #1 \end{pmatrix}}

\def\UCBt{{\tilde{\UCB}}}

\def\BCBov{{\overline{\BCB}}}

\def\piov{{\overline{\pi}}}
\def\Xbov{{\overline{\Xb}}}
\def\Ybov{{\overline{\Yb}}}
\def\ups{{\upsilon}}

\def\upsb{{\boldsymbol{\upsilon}}}

\DeclareMathOperator{\proj}{{\mathrm{proj}}}

\def\and{et}

\def\iso{\buildrel \sim\over\to}

\begin{document}
\title{Compactification des vari\'et\'es 
de Deligne-Lusztig}

\author{C. Bonnaf\'e \& R. Rouquier}
\address{\noindent 
C\'edric BONNAF\'E~: Labo. de Math. de Besan\c{c}on (CNRS: UMR 6623), 
Universit\'e de Franche-Comt\'e, 16 Route de Gray, 25030 Besan\c{c}on
Cedex, France} 
\def\emailaddrname{{\it Courriel~}}
\makeatletter
\email{cedric.bonnafe@univ-fcomte.fr}
%\urladdr{www-math.univ-fcomte.fr/pp\_Annu/CBONNAFE/}

\makeatother

\address{\noindent Rapha\"el ROUQUIER~: Mathematical Institute, 
University of Oxford, 24-29 St Giles', Oxford, OX1 3LB, UK}
\email{rouquier@maths.ox.ac.uk}
%\urladdr{???}
 
%\subjclass{According to the 2000 classification:
%Primary ???; Secondary ???}

\date{\today}

\begin{abstract} 
Nous construisons explicitement la normalisation de la compactification de 
Bott-Samelson-Demazure-Hansen des vari\'et\'es de Deligne-Lusztig $\Xb(w)$ dans
leur rev\^etement $\Yb(w)$ et 
retrouvons ainsi un r\'esultat de Deligne-Lusztig \cite[Lemma 9.13]{delu}
sur la monodromie locale autour des diviseurs de la compactification.
\end{abstract}

\maketitle

\pagestyle{myheadings}

\markboth{\sc C. Bonnaf\'e \& R. Rouquier}{\sc Compactification des vari\'et\'es  
de Deligne-Lusztig}

\tableofcontents

\section*{Introduction}

Dans \cite{jordan}, nous avons \'etudi\'e le prolongement de 
certains syst\`emes 
locaux sur les vari\'et\'es de Deligne-Lusztig en vue d'une application 
alg\'ebrique (\'equivalence de Morita donn\'ee par la d\'ecomposition 
de Jordan, conjectur\'ee par Brou\'e). Dans cette \'etude, nous utilisions
un r\'esultat crucial 
de Deligne-Lusztig sur la ramification de ces syst\`emes locaux 
\cite[lemme 9.13]{delu}. Une des motivations du pr\'esent travail 
est de fournir une alternative ``explicite'' au calcul local effectu\'e
dans la preuve de Deligne et Lusztig. 

Plus pr\'ecis\'ement, si $w$ est un \'el\'ement du groupe de Weyl d'un 
groupe r\'eductif connexe $\Gb$ muni d'une isog\'enie $F$ dont 
une puissance est un endomorphisme de Frobenius, il lui est associ\'e 
deux vari\'et\'es de Deligne-Lusztig $\Xb(w)$ et $\Yb(w)$ 
ainsi qu'un morphisme fini \'etale $\Yb(w) \to \Xb(w)$ faisant de 
$\Xb(w)$ un quotient de $\Yb(w)$ par l'action du groupe fini $\Tb^F$ des
points rationnels d'un tore maximal $F$-stable $\Tb$ de $\Gb$
(voir \cite[\S 1]{delu}~: la vari\'et\'e $\Yb(w)$ y est not\'ee 
$\Xbt(\wdo)$). Deligne et Lusztig 
\cite[lemme 9.11]{delu} 
ont construit une compactification lisse $\Xbov(w)$ de $\Xb(w)$ \`a la
Bott-Samelson-Demazure-Hansen. 
Le but principal de cet article est de construire explicitement la 
{\it normalisation} $\Ybov(w)$ de $\Xbov(w)$ dans $\Yb(w)$~:
$$\diagram
\Yb(w) \dto \rto|<\ahook & \Ybov(w) \dto \\
\Xb(w) \rto|<\ahook & \Xbov(w).
\enddiagram$$
Une fois cette construction explicite r\'ealis\'ee, nous en d\'eduisons 
les propri\'et\'es fondamentales 
de $\Ybov(w)$ (voir le th\'eor\`eme \ref{main}) permettant d'en 
d\'eduire une nouvelle preuve du lemme 9.13 de Deligne-Lusztig \cite{delu}
qui d\'etermine la monodromie locale du rev\^etement le long d'une des
composantes de $\Xbov(w)-\Xb(w)$. Ce lemme est un point clef
dans la preuve de Deligne-Lusztig des conjectures de Macdonald associant une
repr\'esentation
irr\'educible de $\Gb^F$ \`a un caract\`ere en position g\'en\'erale
de $\Tb^F$.

\bigskip

\section*{Notations}

\medskip

Tout au long de cet article, nous fixons un groupe r\'eductif 
connexe $\Gb$ d\'efini sur une cl\^oture alg\'ebrique $\FM$ 
du corps fini \`a $p$ \'el\'ements $\gfp$, o\`u $p$ est un nombre premier. 
Nous supposons de plus que $\Gb$ est muni d'une isog\'enie 
$F : \Gb \to \Gb$ dont une puissance est un endomorphisme de Frobenius 
de $\Gb$. 

Fixons un sous-groupe de Borel $F$-stable $\Bb$ de $\Gb$, un tore 
maximal $F$-stable $\Tb$ de $\Bb$ et notons $\Ub$ le radical 
unipotent de $\Bb$. Notons $W=N_\Gb(\Tb)/\Tb$ le groupe de 
Weyl de $\Gb$ relativement \`a $\Tb$, $X(\Tb)$ (resp. $Y(\Tb)$) 
le r\'eseau des caract\`eres (resp. des sous-groupes \`a un param\`etre) 
de $\Tb$, $\Phi$ (resp. $\Phi^\vee$) le syst\`eme de racines 
(resp. coracines) de $\Gb$ relativement \`a $\Tb$, $\D$ 
(resp. $\D^\vee$) la base de $\Phi$ (resp. $\Phi^\vee$) 
associ\'ee \`a $\Bb$ et $\Phi_+$ (resp. $\Phi_+^\vee$) l'unique 
syst\`eme de racines (resp. coracines) positives contenant $\D$ 
(resp. $\D^\vee$). 

Si $\a \in \Phi$, on notera $\a^\vee$ sa coracine associ\'ee, 
$s_\a \in W$ la r\'eflexion par rapport \`a $\a$, $\Ub_\a$ 
le sous-groupe unipotent \`a un param\`etre normalis\'e 
par $\Tb$ associ\'e \`a $\a$, $\Tb_{\a^\vee}$ le sous-tore de $\Tb$ 
image de $\a^\vee$ et $\Gb_\a$ le sous-groupe de $\Gb$ engendr\'e 
par $\Ub_\a$ et $\Ub_{-\a}$. 

Posons $S=\{s_\a~|~\a \in \D\}$ et $\Sba=S \cup \{1\}$. 
Nous noterons $\ell : W \to \NM=\{0,1,2,\dots\}$ la fonction 
longueur relativement \`a $S$. Nous noterons $B$ le groupe de 
tresses associ\'e \`a $(W,S)$, de g\'en\'erateurs $\{\sb_\a~|~\a \in \D\}$. 
Soit $f : B \to W$ le morphisme canonique (i.e. l'unique morphisme 
tel que $f(\sb_\a)=s_\a$ pour tout $\a \in \D$) et soit $\s : W \to B$ 
l'unique application telle que $\s(s_\a)=\sb_\a$ pour tout $\a \in \D$ 
et $\s(vw)=\s(v)\s(w)$ si $\ell(vw)=\ell(v)+\ell(w)$. Cette 
application v\'erifie $f \circ \s = \Id_W$. 

\bigskip

\section{Vari\'et\'es de Deligne-Lusztig}

\medskip

Le lecteur pourra trouver dans \cite{dimirou} les 
r\'esultats g\'en\'eraux sur les vari\'et\'es de Deligne-Lusztig 
que nous utiliserons ici.

\bigskip

\subsection{D\'efinition} 
Si $n \in N_\Gb(\Tb)$ et si $g\Ub$, $h\Ub \in \Gb/\Ub$, nous \'ecrirons 
$g\Ub \longtrait{n} h\Ub$ si $g^{-1}h \in \Ub n\Ub$. Si $w \in W$ et si 
$g\Bb$, $h\Bb \in \Gb/\Bb$, nous \'ecrirons $g\Bb \longtrait{w} h\Bb$ 
si $g^{-1}h \in \Bb w \Bb$.

Si $\nb=(n_1,\dots,n_r)$ 
est une suite d'\'el\'ements de $N_\Gb(\Tb)$ et si $\wb=(w_1,\dots,w_r)$ 
d\'esigne la suite de leurs images respectives dans $W$, on pose
\eqna
\UCB(\nb)&=&\{(g_1\Ub,\dots,g_r\Ub,g_{r+1}\Ub) \in (\Gb/\Ub)^{r+1}~|~
\\
&&\qquad g_1\Ub \longtrait{n_1} g_2\Ub \longtrait{n_2} \cdots 
\longtrait{n_{r-1}} g_r \Ub \longtrait{n_r} g_{r+1}\Ub\}
\endeqna
et
\eqna
\BCB(\wb)&=&\{(g_1\Bb,\dots,g_r\Bb,g_{r+1}\Bb) \in (\Gb/\Bb)^{r+1}~|~
\\
&&\qquad g_1\Bb \longtrait{w_1} g_2\Bb \longtrait{w_2} \cdots 
\longtrait{w_{r-1}} g_r \Bb \longtrait{w_r} g_{r+1}\Bb\}
\endeqna
Si $t \in \Tb$ et $(g_1\Ub,\dots,g_r\Ub,g_{r+1}\Ub) \in \UCB(\nb)$, 
on pose 
$$(g_1\Ub,g_2\Ub,\dots,g_r\Ub,g_{r+1}\Ub) \cdot t = 
(g_1t\Ub,g_2 \lexp{n_1}{t}\Ub,\dots,g_r \lexp{n_{r-1}\cdots n_1}{t}\Ub,
g_{r+1}\lexp{n_r\cdots n_1}{t}\Ub).$$
Il est alors facile de v\'erifier que, si $\gb \in \UCB(\nb)$, 
alors $\gb\cdot t \in \UCB(\nb)$ et cela d\'efinit 
une action \`a droite de $\Tb$ sur $\UCB(\nb)$. 
De plus, le morphisme canonique $\Gb/\Ub \to \Gb/\Bb$ 
induit un morphisme
$$\fonction{\pib_\nb}{\UCB(\nb)}{\BCB(\wb)}{(g_1\Ub,\dots,g_r\Ub)}{
(g_1\Bb,\dots,g_r\Bb)}$$
et ce dernier induit un isomorphisme
\equat\label{iso quotient}
\UCB(\nb)/\Tb \iso \BCB(\wb).
\endequat
Posons maintenant 
$$\fonction{\upsb_\nb}{\UCB(\nb)}{\Gb/\Ub \times \Gb/\Ub}{
(g_1\Ub,\dots,g_{r+1}\Ub)}{(g_1\Ub,g_{r+1}\Ub)}$$
$$\fonction{\betb_\wb}{\BCB(\wb)}{\Gb/\Bb \times \Gb/\Bb}{
(g_1\Bb,\dots,g_{r+1}\Bb)}{(g_1\Bb,g_{r+1}\Bb).}\leqno{\text{et}}$$
Alors le diagramme
$$\xymatrix{
\UCB(\nb)\ar[d]_{\DS{\pib_\nb}}\ar[rr]^-{\DS{\upsb_\nb}}&&
\Gb/\Ub\times\Gb/\Ub\ar[d]\\
\BCB(\wb)\ar[rr]_-{\DS{\betb_\wb}} &&\Gb/\Bb\times\Gb/\Bb
}$$
est commutatif (la fl\`eche verticale de droite \'etant la projection 
canonique). 

Notons $\UCB_F$ (resp. $\BCB_F$) le graphe du morphisme de Frobenius 
$F : \Gb/\Ub \to \Gb/\Ub$ (resp. $F : \Gb/\Bb \to \Gb/\Bb$). Les 
{\it vari\'et\'es de Deligne-Lusztig} associ\'ees \`a $\nb$ et $\wb$ 
sont respectivement d\'efinies par
$$\Yb(\nb) = \upsb_\nb^{-1}(\UCB_F)\quad\text{et}\quad
\Xb(\wb)=\betb_\wb^{-1}(\BCB_F).$$
Notons toujours $\wb : \Tb \to \Tb$ la conjugaison par 
$w_1 \cdots w_r$. Alors le groupe $\Tb^{\wb F}$ agit 
sur $\Yb(\nb)$ (par restriction de l'action de $\Tb$ sur $\UCB(\nb)$) 
et le morphisme canonique $\pi_\nb : \Yb(\nb) \to \Xb(\wb)$ 
obtenu par restriction de $\pib_\nb$ induit un isomorphisme
\equat\label{quotient F}
\Yb(\nb)/\Tb^{\wb F} \iso \Xb(\wb).
\endequat
Pour finir, notons $\ups_\nb : \Yb(\nb) \to \Gb/\Ub$ 
et $\b_\wb : \Xb(\wb) \to \Gb/\Bb$ les premi\`eres projections. 
Alors le diagramme 
$$\xymatrix{
\Yb(\nb)\ar[d]_{\DS{\pi_\nb}}\ar[rr]^{\DS{\ups_\nb}}&&\Gb/\Ub\ar[d]\\
\Xb(\wb)\ar[rr]_{\DS{\b_\wb}} &&\Gb/\Bb
}$$
est commutatif (la fl\`eche verticale de droite \'etant la projection 
canonique). 

\bigskip

\remark{independance}
Prolongeons l'application $\s : W \to B$ aux suites d'\'el\'ements de 
$W$ en posant $\sigb(\wb)=\s(w_1)\cdots \s(w_r)$. Si $\nb'$ est 
une autre suite d'el\'ements de $N_\Gb(\Tb)$ dont la suite des 
images dans $W$ est $\wb'$, et si $\sigb(\wb)=\sigb(\wb')$, alors 
les vari\'et\'es $\BCB(\wb)$ et $\BCB(\wb')$ sont canoniquement 
isomorphes et les $\Tb$-vari\'et\'es 
$\UCB(\nb)$ et $\UCB(\nb')$ sont isomorphes, ces isomorphismes 
rendant le diagramme
$$\xymatrix{
\UCB(\nb) \ar[drr]_{\DS{\upsb_\nb}}
\ar[dd]_{\DS{\pib_\nb}} \ar[rrrr]^{\DS{\sim}} &&&& \UCB(\nb') 
\ar[dd]^{\DS{\pib_{\nb'}}} \ar[dll]^{\DS{\upsb_{\nb'}}}\\
&&\Gb/\Ub\times\Gb/\Ub \xto'[1,0][2,0] && \\
\BCB(\wb) \ar[rrrr]^{\DS{\sim\qquad\qquad}} \ar[drr]_{\DS{\betb_\wb}}&&~&&
\BCB(\wb')\ar[dll]^{\DS{\betb_{\wb'}}}\\
&&\Gb/\Bb\times\Gb/\Bb && \\
}$$
commutatif. 

De plus, $\Tb^{\wb F}=\Tb^{\wb' F}$ et les $\Tb^{\wb F}$-vari\'et\'es 
$\Yb(\nb)$ et $\Yb(\nb')$ 
(resp. les vari\'et\'es $\Xb(\wb)$ et $\Xb(\wb')$) sont 
isomorphes (resp. canoniquement isomorphes), les isomorphismes 
rendant le diagramme 
$$\diagram
\Yb(\nb) \drrto^{\DS{\ups_\nb}}
\ddto_{\DS{\pi_\nb}} \xto[0,4]^{\DS{\sim}} &&&& \Yb(\nb') 
\ddto^{\DS{\pi_{\nb'}}} \dllto_{\DS{\ups_{\nb'}}}\\
&&\Gb/\Ub \xto'[1,0][2,0] && \\
\Xb(\wb) \xto[0,4]^{\DS{\sim\qquad\qquad}} \drrto^{\DS{\b_\wb}}&&~&&
\Xb(\wb')\dllto_{\DS{\b_{\wb'}}}\\
&&\Gb/\Bb && \\
\enddiagram$$
commutatif.\finl

\bigskip

\subsection{Compactification de Bott-Samelson-Demazure-Hansen}
Pour tout $\a \in \D$, on fixe un repr\'esentant $\sdo_\a$ de 
$s_\a$ dans $\Gb_\a$. La remarque \ref{independance} montre que, 
dans le but de construire une compactification des vari\'et\'es 
$\Yb(\nb)$ et $\Xb(\wb)$, il est suffisant de travailler sous 
les hypoth\`eses suivantes~:

\bigskip

\begin{quotation}
\noindent{\bf Hypoth\`ese :} {\it Nous fixons une suite 
$(\a_1,\dots,\a_r)$ d'\'el\'ements de $\D$ et, si $1 \le i \le r$, 
nous posons pour simplifier $s_i=s_{\a_i}$ et $\sdo_i=\sdo_{\a_i}$. 
Nous supposons de plus que $\nb=(\sdo_1,\dots,\sdo_r)$ et 
$\wb=(s_1,\dots,s_r)$. }
\end{quotation}

\bigskip

Si $\xb=(x_1,\dots,x_r)$ et $\yb=(y_1,\dots,y_r)$ sont deux suites 
d'\'el\'ements de $\Sba=S\cup\{1\}$ (de m\^eme longueur), nous \'ecrirons 
$\xb \infspe \yb$ si, pour tout $i \in \{1,2,\dots,r\}$, 
on a $x_i \in \{1,y_i\}$. On pose aussi $\xbdo=(\xdo_1,\dots,\xdo_r)$, 
o\`u nous choisirons toujours $\dot{1}=1$. Par exemple, $\wbdo=\nb$ et, 
pour simplifier les notations, nous noterons $\pib_\xb$, $\pi_\xb$, 
$\upsb_\xb$ et $\ups_\xb$ les applications $\pib_\xbdo$, $\pi_\xbdo$, 
$\upsb_\xbdo$ et $\ups_\xbdo$, et la vari\'et\'e $\Yb(\xbdo)$ sera not\'ee 
$\Yb(\xb)$. Pour finir, on pose 
$I_\xb=\{1 \le i \le r~|~x_i=1\}$ et on d\'efinit, comme dans 
\cite[\S 4.4.2]{jordan}, 
$$Y_{\wb,\xb}=\sum_{i \in I_\xb} \ZM~s_1 \cdots s_{i-1}(\a_i^\vee).$$

Bott-Samelson, Demazure et Hansen ont construit une compactification lisse
$\BCBov(\wb)$ de $\BCB(\wb)$:
\eqna
\BCBov(\wb)&=&\DS{\coprod_{\xb \infspe \wb} \BCB(\xb)} \\
&=& \{(g_1\Bb,\dots,g_{r+1}\Bb) \in (\Gb/\Bb)^{r+1}~|~\forall~1\le i\le r,~
g_i^{-1} g_{i+1} \in \Gb_{\a_i}\Bb\}.
\endeqna
Alors $\BCBov(\wb)$ est lisse, projective, irr\'eductible et contient 
$\BCB(\wb)$ comme sous-vari\'et\'e ouverte. Posons 
$$\fonction{\overline{\betb}_\wb}{\BCBov(\wb)}{\Gb/\Bb\times\Gb/\Bb}{
(g_1\Bb,\dots,g_{r+1}\Bb)}{(g_1\Bb,g_{r+1}\Bb).}$$
Alors $\overline{\betb}_\wb$ prolonge $\betb_\wb$ (et en fait co\"{\i}ncide 
avec $\betb_\xb$ sur $\BCB(\xb)$ pour tout $\xb \infspe \wb$). 
On pose alors, suivant \cite[\S 9.10]{delu},
$$\Xbov(\wb)=\overline{\betb}_\wb^{-1}(\BCB_F).$$
Notons que
\equat\label{union X}
\Xbov(\wb)=\coprod_{\xb \infspe \wb} \Xb(\xb).
\endequat
Alors $\Xbov(\wb)$ est une vari\'et\'e lisse, projective 
et contient $\Xb(\wb)$ comme sous-vari\'et\'e ouverte \cite[lemme 9.11]{delu}.

\bigskip

\subsection{Normalisation} 
Avant de parler de la compactification de $\Yb(\nb)$ et avant d'\'enoncer 
le r\'esultat principal de cet article, nous aurons besoin de quelques 
notations. Tout d'abord, fixons un entier naturel non nul $d$ et 
une puissance $q$ de $p$ tels que, pour tout $t \in \Tb$ et 
pour tout $w \in W$, on ait $(wF)^d(t)=t^q$. On fixe une racine primitive 
$(q-1)$-i\`eme de l'unit\'e $\zeta$ dans $\Gb_m$. On note encore 
$\wb F : Y(\Tb) \to Y(\Tb)$ l'endomorphisme de groupes induits par 
l'endomorphisme $\wb F : \Tb \to \Tb$ et on pose
$$\fonction{N_\wb}{Y(\Tb)}{\Tb^{\wb F}}{\l}{N_{F^d/\wb F}(\l(\z)),}$$
o\`u $N_{F^d/\wb F} : \Tb \to \Tb$, $t \mapsto t\cdot \lexp{\wb F}{t}~\cdots
~\lexp{(\wb F)^{d-1}}{t}$. Rappelons que $N_\wb$ est surjective et 
induit un isomorphisme 
$$Y(\Tb)/(\wb F - 1)(Y(\Tb)) \iso \Tb^{\wb F}.$$

Le morphisme $\pi_\wb : \Yb(\wb) \to \Xb(\wb)$ \'etant 
fini, on peut d\'efinir la {\it normalisation} $\Ybov(\wb)$ 
de $\Xbov(\wb)$ dans $\Yb(\wb)$~: c'est l'unique vari\'et\'e 
normale $\Zb$ contenant $\Yb(\wb)$ comme sous-vari\'et\'e ouverte 
dense et munie d'un morphisme fini $\piov_\wb : \Zb \to \Xbov(\wb)$ 
prolongeant $\pi_\wb$. Le morphisme 
$\piov_\wb : \Ybov(\wb) \to \Xbov(\wb)$ \'etant fini, 
$\Ybov(\wb)$ est une vari\'et\'e projective. Le but de cet 
article est de la construire explicitement et d'en d\'eduire 
les propri\'et\'es suivantes~:

\bigskip

\begin{theo}\label{main}
Avec les notations pr\'ec\'edentes, on a~:
\begin{itemize}
\itemth{a} La vari\'et\'e $\Ybov(\wb)$ est une vari\'et\'e 
projective, normale, rationnellement lisse, 
de lieu singulier contenu dans 
$$\piov_\wb^{-1}\Bigl(\bigcup_{\stackrel{\xb \infspe \wb}{|I_\xb| \ge 2}}
\Xb(\xb)\Bigr).$$

\itemth{b} La vari\'et\'e $\Ybov(\wb)$ est munie d'une action de 
$\Tb^{\wb F}$ prolongeant l'action sur $\Yb(\wb)$ et telle que 
$\piov_\wb$ induit un isomorphisme 
$\Ybov(\wb)/\Tb^{\wb F} \iso \Xbov(\wb)$. 

\itemth{c} Si $\xb \infspe \wb$, le stabilisateur dans $\Tb^{\wb F}$ 
d'un \'el\'ement de $\piov_\wb^{-1}(\Xb(\xb))$ est 
\'egal \`a $N_\wb(Y_{\wb,\xb})$.

\itemth{d} Si $\xb \infspe \wb$, alors il existe un morphisme 
canonique $i_\xb : \Yb(\xb) \to \piov_\wb^{-1}(\Xb(\xb))$ 
rendant le diagramme suivant commutatif
$$\diagram
\Yb(\xb) \rrto^{i_\xb}\ddrrto_{\DS{\pi_\xb}} && 
\piov_\wb^{-1}(\Xb(\xb)) \ddto^{\DS{\piov_\wb}} \\
&&\\
&& \Xb(\xb)
\enddiagram$$
et induisant un isomorphisme $\Yb(\xb)/N_\xb(Y_{\wb,\xb}) \iso 
\piov_\wb^{-1}(\Xb(\xb))$.
\end{itemize}
\end{theo}

\bigskip

Rassemblons les constructions pr\'ec\'edentes dans le diagramme 
commutatif suivant~:
\equat\label{diag}
\diagram
\Yb(\wb) \xto[0,2]^{\DS{i_\wb}}|<\ahook \ddto_{\DS{\pi_\wb}} && 
\Ybov(\wb) \ddto_{\DS{\piov_\wb}} && 
\piov_\wb^{-1}(\Xb(\xb)) \xto[0,-2]|<\bhook \ddto && 
\Yb(\xb) \xto[0,-2]_{\DS{i_\xb}} \ddllto^{\DS{\pi_\xb}} \\
&&&&&&\\
\Xb(\wb) \xto[0,2]|<\ahook && \Xbov(\wb) && \Xb(\xb) \xto[0,-2]|<\bhook
\enddiagram
\endequat

\remark{rappel iso} 
L'\'enonc\'e (c) du th\'eor\`eme pr\'ec\'edent montre que 
$$\Xb(\xb) \simeq \piov_\wb^{-1}(\Xb(\xb))/
\bigl(\Tb^{\wb F}/N_\wb(Y_{\wb,\xb})\bigr)$$
tandis que l'\'enonc\'e (d) montre que 
$$\Xb(\xb) \simeq \piov_\wb^{-1}(\Xb(\xb))/
\bigl(\Tb^{\xb F}/N_\xb(Y_{\wb,\xb})\bigr).$$
Ceci n'est pas une incoh\'erence car
$$\Tb^{\wb F}/N_\wb(Y_{\wb,\xb}) \simeq \Tb^{\xb F}/N_\xb(Y_{\wb,\xb})$$ 
d'apr\`es \cite[proposition 4.4 (4)]{jordan}.\finl

\bigskip

La section suivante est consacr\'ee \`a la d\'emonstration du th\'eor\`eme 
\ref{main}. Avant cela, montrons que 
ce th\'eor\`eme fournit une autre preuve de \cite[lemme 9.13]{delu}. 
Tout d'abord, si $1 \le i \le r$, 
notons $\wb(i)=(s_1,\dots,s_{i-1},1,s_{i+1},\dots,s_r) \infspe \wb$. 
Alors 
$$\Xbov(\wb) \setminus \Xb(\wb) = \bigcup_{i=1}^r \overline{\Xb(\wb(i))}$$
et les $\overline{\Xb(\wb(i))}$ sont des diviseurs 
lisses \`a croisements normaux.

\bigskip

On en d\'eduit alors \cite[Lemma 9.13]{delu}:

\begin{coro}\label{913}
Le $\Tb^{\wb F}$-torseur $\Yb(\wb)$ (au-dessus de $\Xb(\wb)$) 
se ramifie le long de $\Xb(\wb(i))$ de la m\^eme fa\c{c}on que le changement
de base sous $s_1\cdots s_{i-1}(\a_i^\vee) : \Gb_m \to \Tb$ 
du rev\^etement de Lang $\Tb \to \Tb$, $t \mapsto t^{-1}\cdot\lexp{\wb F}{t}$,
se ramifie en $0$. 
\end{coro}

\bigskip

\section{D\'emonstration du th\'eor\`eme \ref{main}}

\medskip

\subsection{Premi\`ere r\'eduction} 
La preuve que nous proposons du th\'eor\`eme \ref{main} passe 
par une construction explicite de $\Ybov(\wb)$. 
Cependant, pour simplifier cette construction, il convient de remarquer 
qu'en raisonnant comme dans \cite[\S 6.2]{jordan}, on peut supposer (et nous 
le ferons) que l'hypoth\`ese suivante est satisfaite~:

\bigskip

\begin{quotation}
\noindent{\bf Hypoth\`ese.} {\it Dor\'enavant, et 
ce jusqu'\`a la fin de \S 2, nous supposerons que le groupe 
d\'eriv\'e de $\Gb$ est simplement connexe.}
\end{quotation}

\bigskip

Notons que ceci implique que $\Gb_\a \simeq \Sb\Lb_2$ et que $\a^\vee$ est 
injective pour toute racine $\a$ (en particulier, $Y(\Tb)/\ZM\a^\vee$ est 
sans torsion).

\bigskip

\subsection{Fonctions bi-invariantes sur $\Gb_\a\Ub$} 
Avant de proc\'eder \`a la construction explicite de $\Ybov(\wb)$, nous 
aurons besoin de quelques r\'esultats pr\'eliminaires sur les fonctions 
r\'eguli\`eres sur $\Gb_\a \Ub$ invariantes par l'action de $\Ub \times \Ub$ 
par translations \`a gauche et \`a droite (ici, $\a$ est une racine simple). Commen\c{c}ons par \'etudier le cas du groupe $\Sb\Lb_2$. 

Notons 
$$\fonction{\ph}{\Sb\Lb_2}{\Ab^1}{\matrice{a & b \\ c & d}}{c}$$
et notons $\Ub_2$ le sous-groupe de $\Sb\Lb_2$ form\'e des 
matrices unipotentes triangulaires sup\'erieures. Il est alors 
facile de v\'erifier que $\ph$ est invariante par l'action 
de $\Ub_2 \times \Ub_2$ sur $\Sb\Lb_2$ (par translations 
\`a gauche et \`a droite). En fait, en notant $\Bb_2$ le groupe 
des matrices triangulaires sup\'erieures de $\Sb\Lb_2$, on a~:

\bigskip

\begin{prop}\label{fourre tout}
Soient $g \in \Sb\Lb_2$, $z \in \Gb_m$, $t=\diag(z,z^{-1})$ 
et $s=\matrice{0 & -1 \\ 1 & 0}$. Alors~:
\begin{itemize}
\itemth{a} $\FM[\Sb\Lb_2]^{\Ub_2 \times \Ub_2}=\FM[\ph]$.

\itemth{b} $\ph(tg)=z^{-1}\ph(g)$ et $\ph(gt)=z \ph(g)$.

\itemth{c} $\ph(t^{-1}g\lexp{s}{t})=\ph(g)$.

\itemth{d} On a $\ph(g)=0$ si et seulement si $g \in \Bb_2$.

\itemth{e} On a $\ph(g)=1$ si et seulement si $g \in \Ub_2 s \Ub_2$.
\end{itemize}
\end{prop}

\bigskip

\begin{proof}
(a) Soit $\psi \in \FM[\Sb\Lb_2]^{\Ub_2 \times \Ub_2}$. 
Il existe un unique polyn\^ome $P \in \FM[T]$ tel que, pour tout 
$c \in \Ab^1$, $\psi\matrice{1 & 0 \\ c & 1}=P(c)$. 
Alors $\psi-P(\ph)$ est une fonction $\Ub_2 \times \Ub_2$-invariante 
sur $\Sb\Lb_2$ et nulle sur $\lexp{s}{\Ub_2}$. Par cons\'equent, 
elle est nulle sur $\Ub_2 \lexp{s}{\Ub_2} \Ub_2$~: or, cet ensemble 
est dense dans $\Sb\Lb_2$, donc $\psi-P(\ph)=0$.

\medskip

(b), (c) (d) et (e) d\'ecoulent de calculs \'el\'ementaires.
\end{proof}

\bigskip

Revenons aux groupes $\Gb_\a \Ub$. 
Fixons une racine simple $\a \in \D$. Choisissons 
un isomorphisme $\aleph_\a : \Sb\Lb_2 \iso \Gb_\a$ de sorte 
que 
$$\aleph_\a(\Ub_2)=\Ub_\a,\quad\aleph_\a(s)=\sdo_\a,\quad 
\text{et}\quad\aleph_\a\matrice{z & 0 \\ 0 & z^{-1}}=\a^\vee(z)$$
pour tout $z \in \Gb_m$. 
Notons $\Ub_\a^*$ le sous-groupe de $\Ub$ engendr\'e par la famille 
$(\Ub_\b)_{\b \in \Phi^+ \setminus \{\a\}}$. On a alors 
$\Gb_\a\Ub = \Gb_\a\Ub_\a^* = \Gb_\a \ltimes \Ub_\a^*$. On note 
$\t_\a : \Gb_\a \ltimes \Ub_\a^* \to \Gb_\a$ la projection 
naturelle. Notons pour finir 
$\ph_\a$ la composition $\ph \circ \aleph_\a^{-1} \circ \t_\a : 
\Gb_\a\Ub \longto \Ab^1$, de sorte que le diagramme 
$$\diagram
\Gb_\a\Ub \rrto^{\DS{\t_\alpha}}\xto[2,4]^{\DS{\ph_\a}} && 
\Gb_\a \rrto^{\DS{\aleph_\a^{-1}}} && \Sb\Lb_2 \ddto^{\DS{\ph}} \\
&&&&\\
&&&&\Ab^1
\enddiagram$$
soit commutatif. 
C'est une fonction r\'eguli\`ere sur $\Gb_\a\Ub$. Notons de plus que 
\equat\label{phi=1}
\ph_\a(\sdo_\a)=1,
\endequat
car $\ph(s)=1$.

\bigskip

\begin{prop}\label{proprietes phi}
Soient $g \in \Gb_\a\Ub$, $u$, $v \in \Ub$, $t \in \Tb$ et $z \in \Gb_m$. 
Alors
\begin{itemize}
\itemth{a} $\ph_\a(ugv)=\ph_\a(g)$.

\itemth{b} $\ph_\a(g \a^\vee(z))=z \ph_\a(g)$ et 
$\ph_\a(\a^\vee(z)g)=z^{-1} \ph_\a(g)$.

\itemth{c} $t^{-1} g~\lexp{s_\a}{t} \in \Gb_\a\Ub$ et 
$\ph_\a(t^{-1} g~ \lexp{s_\a}{t})=\ph(g)$.

\itemth{d} $\ph_\a(g)=0$ si et seulement si $g \in \Bb$ 
(c'est-\`a-dire si et seulement si $g \in \Tb_{\a^\vee} \Ub = \Bb \cap \Gb_\a\Ub$). 

\itemth{e} $\ph_\a(g)=1$ si et seulement si $g \in \Ub \sdo_\a \Ub$.
\end{itemize}
\end{prop}

\bigskip

\begin{proof}
Les assertions (a), (b), (d) et (e) d\'ecoulent facilement 
de la proposition \ref{fourre tout} et du fait que $\Ub_\a^*$ 
est normalis\'e par $\Ub$ et $\Gb_\a$. Seul le (c) n\'ecessite un 
commentaire. Tout d'abord, comme $\Tb$ est engendr\'e par 
$\Ker \a$ et $\Tb_{\a^\vee}$, il suffit de montrer le r\'esultat 
dans les deux cas suivants~: $\a(t)=1$ ou $t=\a^\vee(z)$, $z \in \Gb_m$. 
Le deuxi\`eme cas se traite imm\'ediatement par la proposition 
\ref{fourre tout}. Dans le premier cas, on remarque 
que $t$ commute avec $\Gb_\a$ (et donc $\lexp{s_\a}{t}=t$) et, 
comme il normalise $\Ub_\a^*$, on a $\t_\a(t^{-1}g~\lexp{s_\a}{t})=\t_\a(g)$.
\end{proof}

\bigskip

\subsection{Construction de la vari\'et\'e ${\boldsymbol{\Ybov(\wb)}}$~: 
premi\`ere \'etape}
Posons pour commencer 
$$\UCBt(\wb)=\{(g_1\Ub,\dots,g_{r+1}\Ub) \in (\Gb/\Ub)^{r+1}~|~
\forall~1\le i\le r,~g_i^{-1} g_{i+1} \in \Gb_{\a_i}\Ub\}$$
et notons 
$$\pibt_\wb : \UCBt(\wb) \longto \BCBov(\wb)$$
l'application canonique.
%Soit 
%$$\Sb(\wb)=\{(t_1,\dots,t_{r+1}) \in \Tb^{r+1}~|~\forall~1 \le i \le r,~
%t_i^{-1} t_{i+1} \in \Tb_{\a_i^\vee}\}.$$
%Alors $\Sb(\wb)$ agit librement \`a droite sur $\UCBov^+(\wbdo)$ de 
%la fa\c{c}on naturelle et il est facile de v\'erifier que 
%$\pibov_\wb^+$ induit un isomorphisme
%\equat\label{piov iso}
%\UCBov^+(\wbdo)/\Sb(\wb) \simeq \BCBov(\wb).
%\endequat
La vari\'et\'e $\UCBt(\wb)$ est irr\'eductible, quasi-affine, lisse 
et de dimension $2r+\dim \Gb/\Ub$. 

Nous d\'efinissons
$$\fonction{\ph_\wb}{\UCBt(\wb)}{\Ab^r}{(g_1\Ub,\dots,g_{r+1}\Ub)}{
\bigl(\ph_{\a_1}(g_1^{-1}g_2),\dots,\ph_{\a_r}(g_r^{-1}g_{r+1})\bigr).}$$
D'apr\`es la proposition \ref{proprietes phi} (a), l'application 
$\ph_\wb$ est bien d\'efinie et est un morphisme de vari\'et\'es. 
Fixons maintenant un $r$-uplet d'entiers naturels non nuls 
$\db=(d_1,\dots,d_r)$, notons $f_\db : \Ab^r \to \Ab^r$, 
$(\x_1,\dots,\x_r) \mapsto (\x_1^{d_1},\dots,\x_r^{d_r})$ et posons 
$$\UCBt_\db(\wb)=\{(\gb,\xib) \in \UCBt(\wb) \times \Ab^r~|~
\ph_\wb(\gb)=f_\db(\xib)\}.$$

\bigskip

\begin{prop}\label{lissite}
La vari\'et\'e $\UCBt_\db(\wb)$ est lisse, de dimension $2r+\dim \Gb/\Ub$.
\end{prop}

\bigskip

\begin{proof}
Si $\a \in \D$ et $d \in \NM^*$, posons 
$$\UCB_{\a,d}=\{(g,\xi) \in \Gb_\a\Ub/\Ub\times\Ab^1~|~\ph_\a(g)=\xi^d\}.$$
Les isomorphismes $\aleph_\alpha:\Gb_\a \iso \Sb\Lb_2$ et
$\Sb\Lb_2/\Ub_2\iso\Ab^2-\{(0,0)\}, \left(\begin{matrix}a&b\\c&d\end{matrix}
\right)\mapsto (a,c)$ induisent un isomorphisme 
$\Gb_\a\Ub/\Ub \iso \Ab^2\setminus\{(0,0)\}$ et finalement
$$\UCB_{\a,d} \iso \{(x,y,\xi) \in \Ab^3~|~(x,y) \neq (0,0)\text{ et } 
y = \xi^d\} \simeq \Ab^2 \setminus \{(0,0)\}.$$
En particulier, $\UCB_{\a,d}$ est lisse. 

Soient $\wb_i=(s_1,\dots,s_i)$ et $\db_i=(d_1,\dots,d_i)$.
On dispose d'une suite de morphismes canoniques 
$$\UCBt_\db(\wbdo)=\UCBt_{\db_r}(\wbdo_r) \longto 
\UCBt_{\db_{r-1}}(\wbdo_{r-1}) \longto \cdots \longto 
\UCBt_{\db_r}(\wbdo_r) \longto \Gb/\Ub$$
(consistant \`a chaque \'etape \`a oublier le dernier terme de $\gb$ et $\xib$) 
qui sont des fibrations successives de
fibres successivement isomorphes \`a des vari\'et\'es de 
la forme $\UCB_{\a,d}$, donc lisses. La lissit\'e de $\UCBt_\db(\wbdo)$
s'en d\'eduit. 
\end{proof}

\bigskip

\subsection{Construction de la vari\'et\'e ${\boldsymbol{\Ybov(\wb)}}$~: 
deuxi\`eme \'etape}
Si $1 \le i \le r$, il existe un unique $\l_i \in Y(\Tb)$ et 
un unique $m_i \in \ZM$ v\'erifiant les 
trois propri\'et\'es suivantes~: 
$$\begin{cases}
\l_i - \wb F(\l_i) = m_i~s_1 \cdots s_{i-1} (\a_i^\vee),&\\
m_i > 0, &\\
\text{$Y(\Tb)/\ZM \l_i$ est sans torsion.}
\end{cases}$$
Ceci d\'ecoule de l'injectivit\'e de $\Id_{Y(\Tb)} - \wb F$ et du 
fait que $Y(\Tb)/\ZM\a_i^\vee$ est sans torsion. 

\bigskip

\remark{p prime} 
Les $m_i$ ne sont pas divisibles par $p$ car 
l'\'egalit\'e qui les d\'efinit implique que 
\equat\label{NW}
m_i N_{F^d/\wb F}(s_1\cdots s_{i-1}(\a_i^\vee)) = (q-1) \l_i,
\endequat 
donc 
\equat\label{mi divise}
m_i \text{ divise } q-1,
\endequat
car $Y(\Tb)/\ZM\l_i$ est sans torsion.\finl

\bigskip

Posons alors $\mb=(m_1,\dots,m_r)$ et 
$$\Ybt(\wb)=\{(\gb;\xib) \in \UCBt_\mb(\wb)~|~
\upsb_\wb(\gb) \in \UCB_F\}.$$
En d'autres termes, $\Ybt(\wb)$ est form\'ee des 
\'el\'ements 
$(g_1\Ub,\dots,g_{r+1}\Ub;\x_1,\dots,\x_r) \in (\Gb/\Ub)^{r+1} \times \Ab^r$ 
tels que
\equat\label{conditions}
\left\{
\begin{array}{l}
\forall~i \in \{1,2,\dots,r\},~g_i^{-1} g_{i+1} \in \Gb_{\a_i}\Ub\text{ et }
\ph_{\a_i}(g_i^{-1} g_{i+1}) = \x_i^{m_i}~;\\
g_{r+1}\Ub=F(g_1)\Ub.
\end{array}\right.
\endequat

Rappelons le lemme suivant~:

\bigskip

\begin{lem}[Deligne-Lusztig]\label{dl}
Soient $\Hb$ un sous-groupe ferm\'e $F$-stable de 
$\Gb$, $\Zb$ une $\Gb$-vari\'et\'e lisse et 
$\th : \Zb \to \Gb/\Hb \times \Gb/\Hb$ un morphisme $\Gb$-\'equivariant. 
Alors le graphe de $F : \Gb/\Hb \to \Gb/\Hb$ est transverse \`a $\th$. 
\end{lem}

\begin{proof}
Ce lemme est montr\'e dans \cite[Preuve du lemme 9.11]{delu} dans le cas 
o\`u $\Hb=\Bb$ mais la preuve reste valable mot pour mot dans le cas 
g\'en\'eral. 
\end{proof}

\bigskip

Le morphisme $\UCBt_\mb(\wb)\to \Gb/\Ub\times\Gb/\Ub$, 
$(\gb,\xib) \mapsto \upsb_\wb(\gb)$ est 
$\Gb$-\'equivariant et la vari\'et\'e $\UCBt_\mb(\wb)$ est lisse 
d'apr\`es la proposition \ref{lissite}. On d\'eduit donc du lemme 
\ref{dl} que 
\equat\label{y lisse}
\text{\it $\Ybt(\wb)$ est lisse, purement de dimension $2r$.}
\endequat
Notons
$$\fonction{\pit_\wb}{\Ybt(\wb)}{\Xbov(\wb)}{(\gb,\xib)}{
\pit_\wb(\gb).}$$

\medskip
Nous allons maintenant construire une action \`a droite de 
$\Tb^{\wb F} \times (\Gb_m)^r$ sur 
$\Ybt(\wb)$. Tout d'abord posons, 
pour tout $\zb=(z_1,\dots,z_r) \in (\Gb_m)^r$,
$$\begin{cases}
\g_1(\zb)=\l_1(z_1)\cdots\l_r(z_r) & \\
\g_{i+1}(\zb)=\lexp{s_i}{\g_i}(\zb) \a_i^\vee(z_i^{m_i}), &
\text{pour $i \in \{1,2,\dots,r\}$.}
\end{cases}$$
Alors
\equat\label{f gamma}
F(\g_1(\zb))=\g_{r+1}(\zb).
\endequat
\begin{proof}
En effet, 
$$\g_{r+1}(\zb)=\lexp{\wb^{-1}}{\g_1(\zb)}~
\lexp{s_r \cdots s_2}{\a_1^\vee(z_1^{m_1})}~
\lexp{s_r \cdots s_3}{\a_2^\vee(z_2^{m_2})}~\cdots~
\a_r^\vee(z_r^{m_r})$$
et donc
\begin{multline*}
 \lexp{\wb F}{\g_1(\zb)}^{-1}~\lexp{\wb}{\g_{r+1}(\zb)} 
= \bigl(\l_1(z_1)\cdots \l_r(z_r)\bigr)~ 
\lexp{\wb F}{\bigl(\l_1(z_1)\cdots \l_r(z_r)\bigr)}^{-1}\\
\times\a_1^\vee(z_1^{-m_1})~\lexp{s_1}{\a_2^\vee(z_2^{-m_2})}~
\cdots ~\lexp{s_1\cdots s_{r-1}}{\a_r^\vee(z_r^{-m_r})} = 1,
\end{multline*}
la derni\`ere \'egalit\'e d\'ecoulant de la d\'efinition des $\l_i$. 
\end{proof}

\bigskip

Si $\zb=(z_1,\dots,z_r) \in (\Gb_m)^r$, si $t \in \Tb^{\wb F}$ et 
si $(g_1\Ub,\dots,g_r\Ub,F(g_1)\Ub;\x_1,\dots,\x_r) \in \Ybt(\wb)$, 
on pose
\begin{multline*}
(g_1\Ub,\dots,g_r\Ub,F(g_1)\Ub;\x_1,\dots,\x_r) * \zb= \\
(g_1\g_1(\zb)\Ub,\dots,g_r\g_r(\zb)\Ub,F(g_1)\g_{r+1}(\zb)\Ub;z_1\x_1,\dots,z_r\x_r)
\end{multline*}
et
\begin{multline*}
(g_1\Ub,\dots,g_r\Ub,F(g_1)\Ub;\x_1,\dots,\x_r)*t = \\
(g_1t\Ub,g_2\lexp{s_1}{t}\Ub,\dots,g_r\lexp{s_{r-1}\cdots s_1}{t}
\Ub,F(g_1)\lexp{s_r\cdots s_1}{t}\Ub;\x_1,\dots,\x_r).
\end{multline*}

\bigskip

\begin{prop}\label{action quotient}
Les formules ci-dessus d\'efinissent une action de 
$\Tb^{\wb F} \times (\Gb_m)^r$ sur la vari\'et\'e $\Ybt(\wb)$. 
De plus, le morphisme $\pit_\wb : \Ybt(\wb) \to \Xbov(\wb)$ 
induit un isomorphisme $\Ybt(\wb)/(\Tb^{\wb F} \times (\Gb_m)^r) \iso 
\Xbov(\wb)$.
\end{prop}

\bigskip

\begin{proof}
Montrons que l'on a bien d\'efini une action. 
Soient $(\gb;\xib) \in \Ybt(\wb)$, $t \in \Tb^{\wb F}$ et 
$\zb \in (\Gb_m)^r$. Il suffit de montrer que $(\gb;\xib)*t$ et 
$(\gb;\xib)*\zb$ appartiennent \`a $\Ybt(\wb)$ (les axiomes des 
actions de groupes sont clairement v\'erifi\'es). \'Ecrivons
$$\gb=(g_1\Ub,\dots,g_{r+1}\Ub),\quad\xib=(\x_1,\dots,\x_r)\quad
\text{et}\quad \zb=(z_1,\dots,z_r).$$

Commen\c{c}ons par montrer que $(\gb;\xib)*t \in \Ybt(\wb)$. 
Posons pour simplifier $t_i=\lexp{s_{i-1}\cdots s_1}{t}$. On a 
alors 
$$(g_i t_i)^{-1}(g_{i+1} t_{i+1}) = t_i^{-1} (g_i^{-1}g_{i+1})~
\lexp{s_i}{t_i}$$
et donc, d'apr\`es la proposition \ref{proprietes phi} (c), on a 
$(g_i t_i)^{-1} (g_{i+1}t_{i+1}) \in \Gb_{\a_i}\Ub$ et 
$$\ph_{\a_i}\bigl((g_it_i)^{-1}(g_{i+1}t_{i+1})\bigr)=
\ph_{\a_i}(g_i^{-1}g_{i+1})=\x_i^{m_i}.$$
D'autre part, puisque $t \in \Tb^{\wb F}$, on a 
$F(g_1 t_1)\Ub=g_{r+1}~\lexp{F}{t}\Ub=g_{r+1}~\lexp{\wb^{-1}}{t}\Ub = 
g_{r+1} t_{r+1} \Ub$. Donc $(\gb;\xib)*t \in \Ybt(\wb)$.

\smallskip

Montrons maintenant que $(\gb;\xib) * \zb \in \Ybt(\wb)$. D'une part, on a
$$(g_i\g_i(\zb))^{-1}(g_{i+1}\g_{i+1}(\zb)) = 
\g_i(\zb)^{-1}(g_i^{-1} g_{i+1}) ~\lexp{s_i}{\g_i(\zb)} \a_i^\vee(z_i^{m_i}),$$
donc $(g_i\g_i(\zb))^{-1}(g_{i+1}\g_{i+1}(\zb)) \in \Gb_{\a_i}\Ub$ et son 
image par $\ph_{\a_i}$ est $z_i^{m_i}\ph_{\a_i}(g_i^{-1}g_{i+1})=(z_i\x_i)^{m_i}$ 
(voir la proposition \ref{proprietes phi} (b) et (c)). D'autre part, 
d'apr\`es \ref{f gamma}, on obtient 
$g_{r+1}\g_{r+1}(\zb) \Ub= F(g_1\g_1(\zb)) \Ub$. 
Donc $(\gb;\xib)*\zb \in \Ybt(\wb)$. 

\medskip

Il nous reste \`a montrer la derni\`ere assertion de la proposition. 
Tout d'abord, il est clair que les $(\Tb^{\wb F} \times (\Gb_m)^r)$-orbites 
sont contenues dans les fibres de $\pit_\wb$. R\'eciproquement, 
montrons que les fibres de $\pit_\wb$ sont des orbites. 
Soient donc $(\gb;\xib)$ et $(\gb';\xib')$ deux \'el\'ements de 
$\Ybt(\wb)$ tels que $\pit_\wb(\gb;\xib)=\pit_\wb(\gb';\xib')$. 
\'Ecrivons 
$$\gb=(g_1\Ub,\dots,g_{r+1}\Ub),\quad \xib=(\x_1,\dots,\x_r),$$
$$\gb'=(g_1'\Ub,\dots,g_{r+1}'\Ub)\quad\text{et}\quad
\xib'=(\x_1',\dots,\x_r').$$
%Quitte \`a faire agir $(\Gb_m)^r$, on peut supposer que 
%$\x_i$ et $\x_i'$ sont \'egaux \`a $0$ ou $1$ pour tout $i$. 
Par hypoth\`ese, il existe $t_i \in \Tb$ tels que 
$g_i' \Ub=g_i t_i \Ub$. Mais, puisque $g_i^{\prime -1} g_{i+1}'$ et 
$g_i^{-1}g_{i+1}$ appartiennent \`a $\Gb_{\a_i}\Ub$, cela montre 
que $t_i^{-1}~t_{i+1}$ appartient \`a $\Tb_{\a_i^\vee}$, 
ou encore que $\lexp{s_i}{t_i}^{-1}~t_{i+1}$ appartient \`a $\Tb_{\a_i^\vee}$. 
Soit donc $z_i \in \Gb_m$ tel que $t_{i+1}=\lexp{s_i}{t_i}~\a_i^\vee(z_i^{m_i})$. 
Posons $\zb=(z_1,\dots,z_r)$. Alors, quitte \`a remplacer 
$(\gb;\xib)$ par $(\gb;\xib)*\zb$, et quitte \`a multiplier $z_i$ 
par une racine $m_i$-i\`eme de l'unit\'e, on peut supposer que 
$\xib=\xib'$ et $t_{i+1}=\lexp{s_i}{t_i}$. Mais alors, le fait 
que $g_{r+1}\Ub=F(g_1)\Ub$ et $g_{r+1}'\Ub=F(g_1')\Ub$ impose 
que $t_{r+1}=F(t_1)$, et donc que $t \in \Tb^{\wb F}$. 
Par cons\'equent, $(\gb';\xib')=(\gb;\xib)*t$. 

\medskip

Les vari\'et\'es $\Ybt(\wb)$ et $\Xbov(\wb)$ \'etant lisses et les fibres 
de $\pit_\wb$ \'etant des $(\Tb^{\wb F} \times (\Gb_m)^r)$-orbites, 
il suffit maintenant
de montrer que $\pit_\wb$ est s\'eparable \cite[proposition 6.6]{borel}. 
Pour cela, notons $\Pi(\xib)=\xi_1\cdots \x_r$ si $(\x_1,\dots,\x_r) \in \Ab^r$ et 
posons 
$$\Ybt_\vide(\wb)=\{(\gb;\xib) \in \Ybt(\wb)~|~\Pi(\xib)\neq 0\}.$$
Notons $\tilde{\iota}_\wb : \Yb(\wb) \to \Ybt(\wb)$, $\gb \mapsto (\gb;1,\dots,1)$. 
Le fait que $\tilde{\iota}_\wb(\gb) \in \Ybt(\wb)$ 
d\'ecoule de ce que $\varphi_{\a_i}(s_i)=1$ 
d'apr\`es \ref{phi=1}. 
Alors $\Ybt_o(\wb)$ est un ouvert de $\Ybt(\wb)$ contenant
$\tilde{\iota}_\wb(\Yb(\wb))$ 
et son image par $\pit_\wb$ est $\Xb(\wb)$ (voir la proposition 
\ref{proprietes phi} (d)). Le morphisme $\Yb(\wb) \to \Xb(\wb)$ \'etant 
s\'eparable, il suffit de remarquer que l'application 
$$\fonctio{\Yb(\wb) \times (\Gb_m)^r}{\Ybt_\vide(\wb)}{(\gb,\zb)}{
\tilde{\iota}_\wb(\gb)*\zb}$$
est un isomorphisme de vari\'et\'es~: cela vient du fait que, 
si $\gb \in (\Gb/\Ub)^{r+1}$ v\'erifie $(\gb;1,\dots,1) \in \Ybt(\wb)$, 
alors $\gb \in \Yb(\wb)$ d'apr\`es la proposition \ref{proprietes phi} (e).
\end{proof}

On a donc un diagramme commutatif, o\`u les fl\`eches verticales sont
des morphismes quotients par les actions des groupes indiqu\'es:

$$\xymatrix{
\Ybt_\vide(\wb) \ar@{^{(}->}[rr]^{\text{ouvert}} \ar[d]_{(\Gb_m)^r} && \Ybt(\wb)
\ar[dd]^{\Tb^{\wb F} \times (\Gb_m)^r} \\
\Yb(\wb) \ar[d]_{\Tb^{\wb F}}\\
\Xb(\wb)\ar@{^{(}->}[rr]_{\text{ouvert}} && \Xbov(\wb)
}$$

\bigskip

Si $I$ est une partie de $\{1,2,\dots,r\}$, on pose 
$$\Ybt_I(\wb)=\{(\gb;\x_1,\dots,\x_r) \in \Ybt(\wb)~|~
\forall~i \in \{1,2,\dots,r\},~i \in I \iff \x_i=0\}$$
et
$$H_I=\{\zb=(z_1,\dots,z_r) \in (\Gb_m)^r~|~
\g_1(\zb)=\cdots=\g_{r+1}(\zb)=1~\text{et}~\forall~i \not\in I,~z_i =1\}.$$
Alors $\Ybt_I(\wb)$ est une sous-vari\'et\'e localement ferm\'ee 
de $\Ybt(\wb)$, stable par l'action de $\Tb^{\wb F} \times (\Gb_m)^r$, 
et 
\equat\label{stab}
\text{\it le stabilisateur d'un \'el\'ement de $\Ybt_I(\wb)$ dans $(\Gb_m)^r$ 
est \'egal \`a $H_I$.}
\endequat
D'autre part,
il r\'esulte facilement de la proposition \ref{proprietes phi} (d) que, 
si $\xb \infspe \wb$, alors
\equat\label{image inverse}
\pit_\wb^{-1}(\Xb(\xb))=\Ybt_{I_{\xb}}(\wb).
\endequat
En particulier, on a une partition en sous-vari\'et\'es localement ferm\'ees
\equat\label{union}
\Ybt(\wb)=\coprod_{I \subset \{1,2,\dots,r\}} \Ybt_I(\wb).
\endequat

\bigskip

\begin{prop}\label{HI fini}
Soit $I$ une partie de $\{1,2,\dots,r\}$. Alors~:
\begin{itemize}
\itemth{a} $H_I$ est un groupe fini, contenu dans $H_{\{1,2,\dots,r\}}$.

\itemth{b} Si $|I| \le 1$, alors $H_I=1$.
\end{itemize}
\end{prop}

\bigskip

\begin{proof}
(a) Il est tout d'abord \'evident que $H_I$ est contenu dans 
$H_{\{1,2,\dots,r\}}$. Il suffit donc de montrer que ce dernier 
est fini. Or, si $\zb=(z_1,\dots,z_r) \in H_{\{1,2,\dots,r\}}$, alors 
puisque $\g_i(\zb)=\g_{i+1}(\zb)=1$, il r\'esulte de la d\'efinition des 
$\g_i$ que $\a_i^\vee(z_i^{m_i})=1$. Puisque $\a_i^\vee$ est injectif 
(car $\Gb_{\a_i} \simeq \Sb\Lb_2$), on en d\'eduit que 
$z_i$ est une racine $m_i$-i\`eme de l'unit\'e. D'o\`u le r\'esultat.

\medskip

(b) Si $I=\vide$, alors $H_I=1$ par d\'efinition. Si $I=\{i\}$ et si 
$\zb=(z_1,\dots,z_r) \in H_I$, alors $z_1=\cdots=z_{i-1}=z_{i+1}=\cdots=z_r=1$. 
Mais de plus $\g_1(\zb)=1$, ce qui implique que $\l_i(z_i)=1$. Donc 
$z_i=1$ car, puisque $Y(\Tb)/\ZM\l_i$ est sans torsion, le morphisme 
$\l_i : \Gb_m \to \Tb$ est injectif.
\end{proof}

\bigskip

Posons maintenant 
$$\Ybov(\wb)=\Ybt(\wb)/(\Gb_m)^r,$$
notons $\proj_\wb : \Ybt(\wb) \to \Ybov(\wb)$ la projection 
canonique et notons $\piov_\wb : \Ybov(\wb) \to \Xbov(\wb)$ le morphisme de 
vari\'et\'es induit par $\pit_\wb$. Si $I \subset \{1,2,\dots,r\}$, 
on note $\xib_I$ la fonction caract\'eristique du compl\'ementaire 
de $I$, que l'on voit comme un \'el\'ement de $\Ab^r$. On pose 
aussi $\Ybov_I(\wb) = \Ybt_I(\wb)/(\Gb_m)^r \subset \Ybov(\wb)$. 
On a bien s\^ur
$$\Ybov(\wb)=\coprod_{I \subset \{1,2,\dots,r\}} \Ybov_I(\wb).$$
Soit $\xb \infspe \wb$.
Comme dans la preuve de la proposition \ref{action quotient}, on montre
qu'on a un morphisme bien d\'efini
\begin{align*}
\Yb(\xb)\times (\Gb_m)^r/H_{I_\xb} &\to \Ybt_{I_\xb}(\wb)\\
(g,z)&\mapsto (g,\xib_{I_\xb})\ast z
\end{align*}
et que c'est un isomorphisme. En particulier, $\Ybt_{I_\xb}(\wb)$ est lisse.

On d\'efinit alors
$$\fonction{i_\xb}{\Yb(\xb)}{\Ybov_{I_\xb}(\wb)}{\gb}{\proj_\wb(\gb,\xib_{I_\xb}).}$$
Il est clair que
\equat\label{piov inverse}
\piov_\wb^{-1}(\Xb(\xb)) = \Ybov_{I_\xb}(\wb).
\endequat
Le morphisme canonique
$\Ybt_{I_\xb}(\wb)\to \Ybov_{I_\xb}(\wb)$ est le quotient par l'action
libre de $(\Gb_m)^r/H_{I_\xb}$
et $\Ybt_{I_\xb}(\wb)$ est lisse, donc
$\Ybov_{I_\xb}(\wb)$ est lisse.

\medskip
Nous allons montrer que $\Ybov(\wb)$ est la normalisation de 
$\Xbov(\wb)$ dans $\Yb(\wb)$ et que les \'enonc\'es (a), (b), (c) et (d) du 
th\'eor\`eme \ref{main} sont v\'erifi\'es.

\bigskip

\subsection{Fin de la d\'emonstration}
Dans la preuve de la proposition \ref{action quotient}, 
il a \'et\'e remarqu\'e que l'application  
$$\fonctio{\Yb(\wb) \times (\Gb_m)^r}{\Ybt_\vide(\wb)}{(\gb,\zb)}{
\tilde{\iota}_\wb(\gb)*\zb}$$
est un isomorphisme de vari\'et\'es. Cela montre que 
$i_\wb : \Yb(\wb) \to \Ybov(\wb)$ est une immersion ouverte, d'image 
$\Ybov_\vide(\wb)=\piov_\wb^{-1}(\Xb(\wb))$. On a donc un 
diagramme commutatif
$$\diagram
\Yb(w) \ddto_{\DS{\pi_\wb}} \rrto^{\DS{i_\wb}}|<\ahook && 
\Ybov(w) \ddto^{\DS{\piov_\wb}} \\
&&\\
\Xb(w) \rrto|<\ahook && \Xbov(w).
\enddiagram$$
D'autre part, par construction, $i_\wb$ est $\Tb^{\wb F}$-\'equivariant 
et il r\'esulte de la proposition \ref{action quotient} 
que $\piov_\wb$ induit un isomorphisme de vari\'et\'es
\equat\label{quotient T}
\Ybov(\wb)/\Tb^{\wb F} \stackrel{\sim}{\longrightarrow} \Xbov(\wb).
\endequat
D'autre part, posons $H=H_{\{1,2,\dots,r\}}$. Alors 
$\Ybt(\wb)/H$ est une vari\'et\'e normale et rationnellement lisse 
(car $\Ybt(\wb)$ est lisse et $H$ est fini) et le groupe 
$(\Gb_m)^r/H$ agit librement sur 
$\Ybt(\wb)/H$. Donc 
\equat\label{normale}
\text{\it $\Ybov(\wb)$ est une vari\'et\'e normale et rationnellement lisse}
\endequat
et
\equat\label{fini}
\text{\it le morphisme $\piov_\wb$ est un morphisme fini.}
\endequat
Par cons\'equent, $\Yb(\wb)$ est bien la normalisation 
de $\Xbov(\wb)$ dans $\Yb(\wb)$.

\bigskip

\subsubsection*{Preuve du (a)} 
Puisque $\piov_\wb$ est un morphisme fini (voir \ref{fini}), 
c'est un morphisme projectif. La vari\'et\'e $\Xbo(\wb)$ \'etant 
projective, $\Ybov(\wb)$ est aussi projective.

Le morphisme canonique
$\Ybt(\wb)/H_{\{1,\ldots,r\}}\to \Ybov(\wb)$ est lisse et la description du
lieu singulier de $\Ybov(\wb)$ se ram\`ene donc au cas de la vari\'et\'e
$\Ybt(\wb)/H_{\{1,\ldots,r\}}$. Puisque $H_{\{1,\ldots,r\}}$ agit
librement sur $\coprod_{|I| \le 1} \Ybt_I(\wb)$
(d'apr\`es la proposition \ref{HI fini} (b)), 
on obtient la derni\`ere assertion de (a). 

\bigskip

\noindent{\sc Remarque - } Le lieu de ramification du morphisme
quotient $\Ybt(\wb)\to\Ybt(\wb)/H_{\{1,\ldots,r\}}$ est de codimension $>1$.
La vari\'et\'e $\Ybt(\wb)$ est lisse, donc le th\'eor\`eme de puret\'e
du lieu de ramification \cite[X, Th\'eor\`eme 3.1]{SGA1} montre que le
lieu singulier de $\Ybt(\wb)/H_{\{1,\ldots,r\}}$ est l'image de
$\coprod_{|H_I| > 1} \Ybt_I(\wb)$. 

Il existe des exemples o\`u la vari\'et\'e normale 
$\Ybov(\wb)$ n'est pas lisse~: 
si $\Gb=\Gb\Lb_3$, si $F$ est l'endomorphisme de Frobenius 
d\'eploy\'e standard sur le corps fini $\FM_{\! q_0}$, et si $s$ et $t$ 
sont les deux r\'eflexions simples, alors la vari\'et\'e $\Ybov(s,t)$ 
n'est pas lisse car le groupe fini $H_{1,2}$ est cyclique d'ordre 
$1+q_0+q_0^2$.\finl

\bigskip

\subsubsection*{Preuve du (b)}
Cela a \'et\'e d\'emontr\'e dans \ref{quotient T}.

\bigskip

\subsubsection*{Preuve du (c)} 
Soit $t \in \Tb^{\wb F}$ et soit $\xb \infspe \wb$. Alors $t$ stabilise 
un \'el\'ement de $\piov_\wb^{-1}(\Xb(\xb))$ si et seulement si il existe 
$\zb=(z_1,\dots,z_r) \in (\Gb_m)^r$ tel que, pour tout $i \in \{1,2,\dots,r\}$, 
on ait
$$\g_i(\zb)=\lexp{s_{i-1}\cdots s_1}{t} \textrm{ pour tout }i
\textrm{ et }z_i=1 \textrm{ pour }i{\not\in}I_\xb.$$
Si tel est le cas, on a 
$z_i^{m_i}=1$ pour tout $i$ et $t=\l_1(z_1)\cdots\l_r(z_r)$. 
D'apr\`es \ref{mi divise}, il existe donc $e_i \in \ZM$ tel que 
$z_i = (\z^{(q-1)/m_i})^{e_i}$, donc 
$$t=N_\wb\bigl(\sum_{i \in I_\xb} e_i~s_1\cdots s_{i-1}(\a_i^\vee)\bigr)$$
d'apr\`es \ref{NW}. 
Donc $t \in N_\wb(Y_{\wb,\xb})$. La r\'eciproque se montre de fa\c{c}on analogue.

\bigskip

\subsubsection*{Preuve du (d)} 
Soit $\xb \infspe \wb$. On a construit un morphisme canonique 
$i_\xb : \Yb(\xb) \to \piov_\wb^{-1}(\Xb(\xb))$. 
Tout d'abord, la surjectivit\'e de $i_\xb$ 
r\'esulte du fait suivant~: si $(\gb;\xib) \in \Ybt(\wb)$ et 
si $\xb \infspe \wb$, alors 
$\xib=\xib_{I_\xb}$ si et seulement si $\gb \in \Yb(\xb)$ 
(voir la proposition \ref{proprietes phi} (e)). 
De plus, puisque 
$\piov_\wb^\xb \circ i_\xb = \pi_\xb$ (o\`u $\piov_\wb^\xb$ d\'esigne 
la restriction de $\piov_\wb$ \`a $\piov_\wb^{-1}(\Xb(\xb)) \to \Xb(\xb)$), 
le morphisme $i_\xb$ est s\'eparable.
Les vari\'et\'es $\Yb(\xb)$ et
$\piov_\wb^{-1}(\Xb(\xb))$ \'etant lisses,
il suffit de montrer 
que les fibres de $i_\xb$ sont les $N_\xb(Y_{\wb,\xb})$-orbites. 

\'Ecrivons $\xb=(x_1,\dots,x_r)$, o\`u $x_i \in \{1,s_i\}$ pour tout $i$. 
Comme $\piov_\wb^\xb \circ i_\xb = \pi_\xb$, les fibres de
$(\piov_\wb^{\xb})^{-1}$ 
sont contenues dans des $\Tb^{\xb F}$-orbites. Soit $t \in \Tb^{\xb F}$ 
et soit $\gb \in \Yb(\xb)$. Alors $i_\xb(\gb \cdot t)=i_\xb(\gb)$ 
si et seulement si il existe $\zb=(z_1,\dots,z_r)$ dans $(\Gb_m)^r$ 
tel que
$$\left\{
\begin{array}{l}
\forall~1 \le i \le r+1,~\lexp{x_{i-1}\cdots x_1}{t}=\g_i(\zb),\\
\forall~i \not\in I_\xb,~z_i=1.
\end{array}\right.\leqno{(*)}$$
Posons, comme dans \cite[\S 4.4.3]{jordan}, 
$$\Sb_{\wb,\xb}=\{(a_1,\dots,a_{r+1}) \in \Tb^{r+1}~|~
a_{r+1}=F(a_1),$$
$$\forall~i {\not\in} I_\xb,~a_{i+1}=\lexp{x_i}{a_i}~ \text{ et }
\forall~i \in I_\xb,~
a_i^{-1}~a_{i+1} \in \Tb_{\a_i^\vee}\}.$$
Alors
$$\Sb_{\wb,\xb}=\{(a_1,\dots,a_{r+1}) \in \Tb^{r+1}~|~
a_{r+1}=F(a_1),$$
$$\forall~i {\not\in} I_\xb,~a_{i+1}=\lexp{x_i}{a_i}~ \text{ et }
\forall~i \in I_\xb,~
\lexp{x_i}{a_i}^{-1}~a_{i+1} \in \Tb_{\a_i^\vee}\}$$
et donc l'application 
$$\fonction{\mu_\xb}{\Tb^{\xb F}}{\Sb_{\wb,\xb}}{t}{(t,\lexp{x_1}{t},\dots,
\lexp{x_{r-1}\cdots x_1}{t},\lexp{x_r\cdots x_1}{t})}$$
est bien d\'efinie (et est un morphisme de groupes injectif). 
D'autre part, si $\Rb_\xb$ est l'ensemble des 
$(z_1,\dots,z_r) \in (\Gb_m)^r$ tels que $z_i=1$ si $i \not\in I_\xb$, 
alors l'application 
$$\fonction{\gamb_\xb}{\Rb_\xb}{\Sb_{\wb,\xb}}{\zb}{
(\g_1(\zb),\dots,\g_{r+1}(\zb))}$$
est un morphisme de groupes qui est bien d\'efini et dont il est facile 
de v\'erifier que le noyau est fini (c'est $H_{I_\xb}$). 
Donc son image est de dimension $|I_\xb|$, ce qui est aussi 
la dimension de $\Sb_{\wb,\xb}$ (voir \cite[page 17]{jordan}). 
Donc 
$$\Sb_{\wb,\xb}^\circ={\mathrm{Im}}(\gamb_\xb).$$
Mais la condition $(*)$ est \'equivalente \`a dire que 
$\mu_\xb(t) \in {\mathrm{Im}}(\gamb_\xb)$. Le r\'esultat 
d\'ecoule alors de \cite[proposition 4.11 (4)]{jordan}.
La preuve du th\'eor\`eme \ref{main} est compl\`ete.

\end{document}